\documentclass{article}
\usepackage{amsfonts}
\usepackage{amsbsy}
\usepackage{amsmath}
\usepackage{amssymb}
\usepackage{amscd}
\usepackage{amsthm}
\usepackage{epsf}
\usepackage{amssymb}
\usepackage{latexsym}

\theoremstyle{plain}
\newtheorem{theorem}{Theorem}[section]
\newtheorem{proposition}[theorem]{Proposition}
\newtheorem{lemma}[theorem]{Lemma}
\newtheorem{corollary}[theorem]{Corollary}
\theoremstyle{definition}

\theoremstyle{remark}

\newtheorem*{rem}{Remark}

\newcommand{\vF} {{\varPhi}}
\newcommand{\om} {{\omega}}
\newcommand{\la} {{\lambda}}
\newcommand{\ve} {{\varepsilon}}
\newcommand{\vf} {{\varphi}}
\newcommand{\De} {{\varDelta}}
\numberwithin{equation}{section}

\begin{document}
\title{On the Divided Power Algebra and the Symplectic Group in Characteristic 2}
\author{Mihalis Maliakas$^*$}
\date{}
\maketitle
\begin{abstract}
Let $V$ be an even dimensional vector space over a field $K$ of
characteristic 2 equipped with a non-degenerate alternating
bilinear form $f$. The divided power algebra $DV$ is considered as
a complex with differential defined from $f$. We examine the
cohomology modules as representations of the corresponding
symplectic group.
\end{abstract}
\footnotetext{

\hspace{0.25cm}Key Words: divided power algebra, symplectic group

\hspace{0.25cm}
* Research supported by ELKE of the University of
Athens}
\section{Introduction}
Let $V$ be a vector space of dimension $2m$ over a field $K$
equipped with a non-degenerate alternating bilinear form $f$. Let
$G=Sp(V)$ be the corresponding symplectic group and
$\bigwedge\limits V$ the exterior algebra of $V$. There are
contraction maps $\theta_k : \bigwedge\limits^k V\to
\bigwedge\limits^{k-2}V$ defined from $f$ that are $G$-equivariant
for the natural action of $G$ on $\bigwedge\limits^k V$. By
extending in the usual way, one obtains a contraction map $\theta$
defined on $\bigwedge\limits V$ which is $G$-equivariant. Assume
now and for the rest of the paper that the characteristic of $K$
is 2. In this case we have $\theta^2 =0$ so that $\bigwedge\limits
V$ becomes a complex. It was shown by R.~Gow in \cite{2} that the
complex $\bigwedge\limits V$ is exact in all degrees except one
and the unique non-vanishing homology module is an irreducible
$G$-module. Moreover, if $K$ is algebraically closed, this module
was identified as the spin module for $G$, that is, an irreducible
rational representation of $G$ of highest weight $\om_m$, where
$\om_1 ,\ldots, \om_m$ are the fundamental weights of $G$ (in the
numbering such that $\om_1$ corresponds to the natural module of
$G$).

The purpose of this note is to examine the case of the divided
power algebra $DV$ of $V$. (We refer to \cite[A2.4]{1} for $DV$).
One reason the divided power algebra arises naturally in the
representation theory of $G$ is that the degree $k$ component $D_k
V$ of $DV=\underset{k\ge 0}{\oplus}D_k V$ is a Weyl module of $G$.
In order to describe our main results we need some notation.

Let $\{ x_1 ,y_1 ,\ldots, x_m ,y_m \}$ be a basis of $V$ such that
$f(x_i ,y_j )=\delta_{ij}$, $f(x_i ,x_j )=f(y_i ,y_j)=0$ for $1\le
i,j\le m$. Since the characteristic of $K$ is 2, there is a well
defined $G$-equivariant map $\bigwedge\limits^2 V\to D_2 V$,
$u\wedge v \mapsto uv$, $u,v \in V$. Let $\om$ be the image of the
element $x_1 \wedge y_1 +\cdots +x_m \wedge y_m$ under this map
and let
\[
    \partial_k :D_k V\to D_{k+2}V
\]
be the map given by multiplication by $\om$. We define a map
$\partial :DV\to DV$ whose restriction to $D_k V$ is $\partial_k$.
It easily follows that $\partial$ is $G$-equivariant for the
natural action of $G$ on $DV$. Since the characteristic of $K$ is
2, we have $\partial^2 =0$. We regard $DV$ as a complex is which
$\partial$ has degree 2 and the cochain module in degree $k$ is
$D_k V$.

The main result of this paper concerns the determination of the
cohomology modules $H^{k}(DV)$. If $k<m$ or $k\not\equiv m \mod 2$
then $H^{k}(DV)=0$ (Corollary 2.3). If $k\ge 0$ and $K$ is
algebraically closed, then there is an isomorphism of $G$-modules,
\[
    H^{m+2k}(DV)\simeq \Delta (k\om_1 )^{(1)}\otimes L(\om_m )
\]
where $\Delta (k\om_1 )$ is the Weyl module of $G$ of highest
weight $k\om_1$, $L(\om_m )$ is the irreducible rational
representation of $G$ of highest weight $\om_m$ and $\Delta
(k\om_1 )^{(1)}$ denotes the first Frobenius twist of $\Delta
(k\om_1 )$ (Theorem 3.2).

The first and second non-vanishing cohomology modules are $H^m
(DV)$ and $H^{m+2}(DV)$. We have that $H^m (DV)$ is an irreducible
$G$-module of highest weight $\om_m$. Thus $H^m (DV)$ is a spin
module for $G$ and the situation here is similar to \cite{2}. We
show that $H^{m+2}(DV)$ is also an irreducible $G$-module and its
highest weight is $2\om_1 +\om_m$ (Corollary 3.4). However, it is
not true in general that the nonzero $H^k (DV)$ are irreducible.
\section{The Complex $DV$}

In this section we will obtain a first description of the
cohomology modules $H^k(DV)=ker\partial_{k}/Im\partial_{k-2}$ of
$DV$.

If $v\in V$ and $a$ is a non-negative integer, we will use the
notation $v^{(a)}$ for the $a$-th divided power of $v$. If $a=0$,
we regard $v^{(0)}=1$. If $\{v_1,...,v_n\}$ is a basis of $V$,
then a basis of $D_kV$ is
$\{v_1^{(a_1)}...v_n^{(a_n)}|a_1+...+a_n=k\}$. There are certain
identities in $DV$ that we will use frequently such as
\[
    v^{(a)}v^{(b)}=\binom{a+b}{a}v^{(a+b)}, \quad
    (tu)^{(a)}=t^au^{(a)},\quad
    (u+v)^{(a)}=\sum^{a}_{i=0}u^{(i)}u^{(a-i)},
\]
where $u,v\in V$ and $t\in K$.

We will need the following simple Lemma in the proof of Theorem
\ref{th2.2} below and elsewhere.
\begin{lemma}\label{le2.1}
Suppose $m=1$. Then the cohomology of the complex $DV$ is given by
\[
    H^k (DV)=\left\{ \begin{array}{cll}
    0, &k& \text{even}\\
    D_k V, &k& \text{odd}.
    \end{array}\right.
\]
\end{lemma}
\begin{proof}
From the definition of the complex $DV$ it follows that
\[
    DV=D^0 V\oplus D^1 V
\]
where $D^0 V$ (respectively, $D^1 V$) is the subcomplex of $DV$
consisting of the cochains of $DV$ that have even (respectively,
odd) degrees, namely
\begin{align*}
    &D^0 V: 0\longrightarrow K \overset{\partial_0}{\longrightarrow}
    D_2 V\overset{\partial_2}{\longrightarrow} D_4 V
    \overset{\partial_4}{\longrightarrow} \cdots\\
    &D^1 V: 0\longrightarrow V \overset{\partial_1}{\longrightarrow}
    D_3 V \overset{\partial_3}{\longrightarrow} D_5 V
    \overset{\partial_5}{\longrightarrow}\cdots
\end{align*}
Hence
\[
    H^k (DV)=\left\{ \begin{array}{lll}
    H^k (D^0 V), & k& \text{even}\\
    H^k (D^1 V), & k& \text{odd}
    \end{array}\right.
\]
Let $m=1$. Suppose $x=x^{(a)}_{1}y^{(b)}_{1}\in D_k V$, $a+b=k$.
Then
\begin{align}
    \partial_k (x)=(a+1)(b+1)x^{(a+1)}_{1}y^{(b+1)}_{1}\tag{1}
\end{align}
\indent Let $k$ be odd. Then one of $a,b$ is odd and since the
characteristic of $K$ is 2, we see from (1) that $\partial_k
(x)=0$ for all $x\in D_k V$. Thus $H^k (DV)=H^k (D^1 V)=D_k V$.

Let $k$ be even. Then $H^k (DV)=H^k (D^0 V)$. We claim that $D^0
V$ is exact. Indeed, it is obvious that $\partial_0$ is injective.
Let  $k\geq2$. It follows from (1) that $\ker
\partial_k$ is spanned by the elements $x^{(a)}_{1}y^{(b)}_{1}$,
$a+b =k$, such that both $a$ and $b$ are odd. For such an element
we have
$\partial_{k-2}(x^{(a-1)}_{1}y^{(b-1)}_{1})=abx^{(a)}_{1}y^{(b)}_{1}=x^{(a)}_{1}y^{(b)}_{1}$
so that ${\rm Im}\,\partial_{k-2}=\ker \partial_k$. Thus $H^k (D^0
V)=0$.
\end{proof}

For each $i=1,\ldots, m$, let $V_i$ be the subspace of $V$ with
basis $\{ x_i ,y_i \}$. Then the $V_i$ are non-degenerate and we
consider the group $H=Sp (V_1 )\times \cdots \times Sp (V_m )$
embedded in $G=Sp (V)$ in the usual manner. Each $G$-module may be
considered an $H$-module upon restriction. In particular the
cohomology modules $H^k (DV)$ of $DV$ are $H$-modules.

\begin{theorem}\label{th2.2}
As $H$-modules we have
\[
    H^k (DV)\simeq \bigoplus D_{a_1}V_1 \otimes \cdots \otimes
    D_{a_m}V_m
\]
where the sum ranges over all $m$-tuples $(a_1 ,\ldots, a_m )$ of
odd positive integers $a_1 ,\ldots, a_m$ such that $a_1 +\cdots +a_m
=k$.
\end{theorem}
\begin{proof}
We use induction on $m$, the case $m=1$ owing to Lemma \ref{le2.1}.

Assume $m\ge 2$ and let $U=V_2 \oplus \cdots \oplus V_m$ so that
$V=V_1 \oplus U$. We have an isomorphism
\[
    D_k V\simeq \bigoplus_{a+b=k}D_a V_1 \otimes D_b U
\]
of $GL(V_1 )\times GL(U)$ modules, where $GL(V)$ denotes the
general linear group of $V$. The injective maps $D_a V_1 \otimes
D_b U\to D_k V$, $a+b=k$, are given by
\[
   D_a V_1 \otimes D_b U\ni z_1 \otimes z_2 \mapsto z_1 z_2 \in D_k V.
\]
Under these
identifications, the differential of the complex $DV$ restricted
to $D_a V_1 \otimes D_b U$ looks like
\[
    \begin{array}{ll}
    D_a V_1 \otimes D_b U & \longrightarrow D_{a+2}V_1 \otimes D_b U\\
    &\;\searrow \quad \quad \quad \;\oplus
    \\
    & \quad \quad D_a V_1 \otimes D_{b+2}U
    \end{array}
\]
where the horizontal map is multiplication by $x_1 y_1$ on the first
factor and the identity on the second and the diagonal map is the
identity on the first factor and multiplication by $x_2 y_2 +\cdots
+x_m y_m$ on the second. It follows that the complex $DV$ is
isomorphic to the tensor product $DV_1 \otimes DU$ of the complexes
$DV_1$ and $DU$ and moreover this isomorphism is $Sp(V_1 )\times
Sp(U)$-equivariant. Thus
\begin{align}
    H^k (DV)\simeq \bigoplus_{a+b=k}H^a (DV_1 )\otimes H^b
    (DU)\tag{2}
\end{align}
as $Sp (V_1 )\times Sp(U)$-modules. From Lemma \ref{le2.1}, $H^a
(DV_1 )=0$ unless $a$ is odd in which case $H^a (DV_1 )=D_a V_1$.
Using this and the induction hypothesis for $H^b (DU)$, the desired
result follows.
\end{proof}

An immediate corollary is the following result.
\begin{corollary}\label{co2.3}
$H^k (DV)=0$ if and only if $k<m$ or $k\not\equiv m \mod 2$.
\end{corollary}

The next result provides a basis for $H^{m+2k}(DV)$, $k\ge 0$,
that will be used in Section 3.
\begin{corollary}\label{re2.4}
A basis for $H^{m+2k}(DV)$, $k\ge 0$, consists of the elements
\begin{align*}
    \ x^{(b_1 )}_{1}y^{(c_1 )}_{1}\cdots
    x^{(b_m )}_{m}y^{(c_m )}_{m}+{\rm Im}\,\partial_{m+2k-2}\
\end{align*}
where $\sum^{m}_{i=1}(b_i +c_i )=
    m+2k$ and $\ b_i +c_i$ is odd for all $\ i$
\end{corollary}
\begin{proof}
The isomorphism in the statement of Theorem \ref{th2.2} yields
injective maps $D_{a_1}V_1 \otimes \cdots \otimes D_{a_m}V_m \to
H^{m+2k}(DV)$, where each $a_i$ is odd and $a_1 +\cdots +a_m =m+2k$.
It follows from the proof of the Theorem, that these are induced
from the injective maps $D_{a_1}V_1 \otimes \cdots \otimes
D_{a_m}V_m \ni z_1 \otimes \cdots \otimes z_m \mapsto z_1 \cdots z_m
\in D_{m+2k} V$, that is, they are given by $D_{a_1}V_1 \otimes
\cdots \otimes D_{a_m}V_m \ni z_1 \otimes \cdots \otimes z_m \mapsto
z_1 \cdots z_m +{\rm Im}\,\partial_{m+2k-2}\in H^{m+2k} (DV)$. A
basis for $D_{a_1}V_1 \otimes \cdots \otimes D_{a_m}V_m$ is $\{
x^{(b_1 )}_{1}y^{(c_1 )}_{1} \otimes \cdots \otimes x^{(b_m
)}_{m}y^{(c_m )}_{m} \mid b_i +c_i =a_i$ for all $i\}$. The union of
the images of these bases under the previous injections yields the
basis of $H^{m+2k}(DV)$ in the statement of the Corollary.
\end{proof}

We determine the dimensions of the nonzero cohomology modules of
$DV$.
\begin{corollary}\label{pro2.5}
Let $k\ge 0$. Then
\[
    \dim H^{m+2k}(DV)=2^m \binom{2m+k-1}{k}.
\]
\end{corollary}
\begin{proof}
We shall make use of the identity
\begin{align}
    \binom{n}{k} =\sum^{k}_{j=0}(j+1) \binom{n-2-j}{k-j},\tag{3}
\end{align}
which may be proved easily by induction on $n$.

Let $f(m,k)=\dim H^{m+2k}(DV)$. We will show by induction on $m$
that
\begin{align}
    f(m,k)=2^m \binom{2m+k-1}{k}.\tag{4}
\end{align}

Let $m=1$. Then from Lemma \ref{le2.1} if follows that
$f(1,k)=\dim D_{2k+1}V=2k+2$.

Let $m\ge 2$. From the isomorphism (2) in the proof of Theorem
\ref{th2.2} and from Lemma \ref{le2.1} we have
\[
    H^{m+2k}(DV)\simeq \bigoplus_{a \ \text{odd}} D_a V_1 \otimes
    H^{m+2k-a} (DU).
\]
By Corollary 2.3, we have $H^{m+2k-a} (DU)=0$, if $a>2k+1$. Hence
in the above sum, $a$ ranges over the odd integers 1,3,...,$2k+1$.
Thus
\[
    f(m,k)=\sum_{a \ \text{odd}}(a+1)f\Big(
    m-1,k-\frac{a-1}{2}\Big),
\]
where $a$ ranges as before. By induction $f(m-1,k-\frac{a-1}{2}) =
2^{m-1}
    \binom{2m+k-\frac{a-1}{2}-3}{k-\frac{a-1}{2}}$
and upon substitution we obtain
\begin{align}
    f( m,k) = 2^{m} \sum_{a \ \text{odd}} \frac{a+1}{2}
    \binom{2m+k-\frac{a-1}{2}-3}{k-\frac{a-1}{2}}.\tag{5}
\end{align}
Using (3), we see that (5) yields (4).
\end{proof}
\section{Main result}
We will consider here the non-vanishing cohomology modules of
$DV$, that is $H^{m+2k}(DV)$, $k\ge 0$, as $G$-modules. In this
section we assume that $K$ is algebraically closed.

First we introduce some notation. Having fixed the ordered basis
$x_1 ,\ldots, x_m ,\linebreak y_m ,\ldots, y_1$ of $V$ we identify
$G= Sp(V)$ with $Sp(2m,K)=\{ A\in GL(2m,K)\mid A^t JA=J\}$ where
$J$ is the matrix of the bilinear form $f$ of the Introduction
with respect to the above basis. A maximal torus in $Sp(2m, K)$ is
$T=\{ diag (t_1 ,\ldots, t_m , t^{-1}_{m},\ldots, t^{-1}_1 )\mid
t_i \neq 0$ for all $i\}$. For each $i=1,\ldots, 2m$ let $\ve_i
:T\to K-\{ 0\}$ be the map the sends a matrix in $T$ to its
element in position $(i,i)$. Then in the weight lattice of $T$ we
have the relations
\[
    \ve_i +\ve_{2m +1-i}=0,
\]
where $1\le i\le 2m$. The set of roots of $Sp (2m,K)$ is $\vF =\{
\pm \ve_i \pm \ve_j ,\pm 2\ve_k \mid 1\le i< j \le m,1\le k\le
m\}$. Let $\varPi = \{ \ve_i -\ve_{i+1}\mid 1\le i<m\} \cup \{
2\ve_m \}$. Then $\varPi$ is a set of simple roots in $\vF$ and
the corresponding positive roots are $\vF^+ = \{ \ve_i +\ve_j
,\ve_i -\ve_j ,2\ve_k \mid 1\le i<j\le m, 1\le k\le m\}$. The
fundamental dominant weights are $\om_1 ,\ldots, \om_m$ where
$\om_i =\ve_1 +\cdots +\ve_i$.

The irreducible rational representations of $Sp(2m,k)$ are
parametrized by dominant weights [3, II.2]. Let $L(\la)$ denote
the (unique up to isomorphism) irreducible module of $Sp(2m,k)$ of
highest weight $\la$, where $\la$ is a dominant weight.

The Weyl module of $G$ of highest weight $\la$, where $\la$ is a
dominant weight, will be denoted by $\Delta (\la)$. It is well
known that for any integer $k\geq0$ we have $\Delta (k\om_1
)\simeq D_k V$ \cite[II 2.17]{3}.
\begin{proposition}\label{pro3.1}
We have $H^m (DV)\simeq L(\om_m )$ as $G$-modules.
\end{proposition}
\begin{proof}
From Theorem \ref{th2.2}, there is an isomorphism of $H$-modules
$H^m (DV)\simeq V_1 \otimes \cdots \otimes V_m$. Since $V_1
\otimes \cdots \otimes V_m$ is irreducible as an $H$-module, it
follows that $H^m (DV)$ is irreducible as a $G$-module. By
Corollary 2.5 a basis of $H^m (DV)$ consists of the elements $v_1
\cdots v_m +{\rm Im}\,\partial_{m-2}$, where for each $i$ we have
$v_i =x_i$ or $y_i$. These are weight vectors for the action of
$G$ and we see that the weights of $H^m (DV)$ are $\pm \ve_1 \pm
\cdots \pm \ve_m$. It is easy to verify that among these, $\ve_1
+\cdots +\ve_m$ is the highest weight. Thus $H^m (DV)\simeq
L(\om_m )$.
\end{proof}

If $M$ is a $G$-module, we use the notation $M^{(1)}$ for the first
Frobenius twist of $M$. The main result of this paper is the
following.
\begin{theorem}\label{th3.2}
Let $k\ge 0$. There is an isomorphism of $G$-modules
\[
    H^{m+2k}(DV)\simeq \Delta (k\om_1 )^{(1)}\otimes L(\om_m ).
\]
\end{theorem}
For the proof we will need the following Lemma.
\begin{lemma}\label{le3.3}
Suppose $V$ be a vector space of dimension $n$ over a field $K$ of
characteristic 2 and $\{ v_1 ,\ldots, v_n \}$ a basis of $V$. Let
$k\ge 0$ and define $N$ to be the subspace of $D_{2k}V$ spanned by
the monomials $v^{(a_1 )}_{1}\cdots v^{(a_n )}_{n}$, $a_1 +\cdots
+ a_n =2k$, where at least one of the $a_i$ is odd. Then $N$ is a
$SL(n,K)$-submodule of $D_{2k}V$ and the map $\vf :(D_k
V)^{(1)}\to D_{2k}V/N$, $v^{(b_1 )}_{1}\cdots v^{(b_n )}_{n}
\mapsto v^{(2b_1 )}_{1} \cdots v^{(2b_n )}_{n} +N$ is an
isomorphism of $SL(n,K)$-modules.
\end{lemma}
\begin{proof}
For $t\in K$, $r,s\in \{ 1,\ldots, n\}$, $r\neq s$, let
$g_{rs}(t)\in SL(n,K)$ be $g_{rs}(t)=I+tE_{rs}$, where $I$ is the
$n\times n$ identity matrix and $E_{rs}$ in the $n\times n$ matrix
with $1$ in position $(r,s)$ and $0$ elsewhere. Since the
$g_{rs}(t)$ generate $SL(n,K)$, in order to prove the first
conclusion of the lemma, it suffices to show that
$g_{rs}(t)N\subseteq N$ for every $t\in K$, $r\neq s$.

We will make use of the fact that $\binom{a+b}{a}$ is even if both
$a,b$ are odd. Let $v=v^{(a_1 )}_{1}\cdots v^{(a_n )}_{n} \in
D_{2k}V$, $a_1 +\cdots + a_n =2k$, such that at least one $a_i$ is
odd. We have
\begin{align}
    g_{rs}(t)v=v^{(a_1 )}_{1}\cdots (v_s +tv_r )^{(a_s )} \cdots
    v^{(a_n )}_{n}\notag \\
    =v^{(a_1 )}_{1} \cdots \bigg( \sum^{a_s}_{i=0}
    t^{a_s -i} v^{(i)}_{s}v^{(a_s -i)}_{r}\bigg)& \cdots v^{(a_n
    )}_{n}\notag.
\end{align}
Thus
\begin{align}
    g_{rs}(t)v=\sum^{a_s}_{i=0}t^{a_s -i}
    \binom{a_r +a_s -i}{a_r} v^{(a_1 )}_{1}\cdots v^{(i)}_{s}
    \cdots v^{(a_r +a_s -i)}_{r} \cdots v^{(a_n )}_{n}.\tag{6}
\end{align}
Consider a monomial $u=v^{(a_1 )}_{1}\cdots v^{(i)}_{s}\cdots
v^{(a_s +a_r -i)}_{r}\cdots v^{(a_n )}_{n}$ in the right hand side
of (6) where all exponents are even. By the assumption on $v$, at
least one of $a_r$, $a_s$ is odd. We conclude that both of $a_r$ and
$a_s -i$ are odd. But then the coefficient $\binom{a_r +a_s
-i}{a_r}$ of $u$ is even and hence zero in $K$. Thus $g_{rs}(t)v\in
N$.

We show now that the map $\vf :(D_k V)^{(1)}\to D_{2k}V/N$ is a
map of $SL(n,K)$-modules. We will use the fact that
$\binom{a}{b}\equiv \binom{2a}{2b} \mod 2$ for all integers $a\ge
b\ge 0$. Let $v=v^{(a_1 )}_{1}\cdots v^{(a_n )}_{n}\in D_k V$. By
a similar computation as in the proof of (6) and taking into
account the Frobenius twist on $D_k V$ we have
\[
    g_{rs}(t)v=\sum^{a_s}_{i=0}t^{2(a_s -i)}
    \binom{a_r +a_s -i}{a_r} v^{(a_1 )}_{1}\cdots v^{(i)}_{s}
    \cdots v^{(a_r +a_s -i)}_{r} \cdots v^{(a_n )}_{n}
\]
so that
\begin{align}
    \vf &(g_{rs}(t)v)=\notag \\
    &=
    \sum^{a_s}_{i=0}t^{2(a_s -i)}
    \binom{a_r +a_s -i}{a_r} v^{(2a_1 )}_{1}\cdots v^{(2i)}_{s}
    \cdots v^{(2(a_r +a_s -i))}_{r} \cdots v^{(2a_n )}_{n}+N.\tag{7}
\end{align}
On the other hand applying (6) to $\vf(v)$ in place of $v$ we have
\begin{align}
    &g_{rs}(t)\vf(v)=\notag \\
    &=
    \sum^{2a_s}_{i=0}t^{2a_s -i}
    \binom{2a_r +2a_s -i}{2a_r} v^{(2a_1 )}_{1}\cdots v^{(i)}_{s}
    \cdots v^{(2a_r +2a_s -i)}_{r} \cdots v^{(2a_n )}_{n}+N\notag\\
    &=\sum^{a_s}_{i=0}t^{2a_s -2i}
    \binom{2a_r +2a_s -2i}{2a_r} v^{(2a_1 )}_{1}\cdots v^{(2i)}_{s}
    \cdots v^{(2a_r +2a_s -2i)}_{r} \cdots v^{(2a_n )}_{n}+N, \
    (8)\notag
  \end{align}
where the last equality follows from the definition of $N$. From
(7), (8) and the fact that $\binom{a}{b}\equiv \binom{2a}{2b}\mod
2$ we see that $\vf (g_{rs}(t)v)=g_{rs}(t)\vf (v)$. It follows
that $\vf$ is a map of $SL(n,k)$-modules.

Finally, it is clear that the map $\vf$ carries a basis of $(D_k
V)^{(1)}$ to a basis of $D_{2k}V/N$ and thus is an isomorphism.
\end{proof}
\begin{proof}[Proof of Theorem \ref{th3.2}]
The $GL(V)$-map $D_{2k}V\otimes D_m V\to D_{m+2k}V$, $u\otimes
v\mapsto uv$, induces a map $D_{2k}V\otimes \ker
\partial_m \to \ker
\partial_{m+2k}$ because if $\partial_m (v)=0$, then $v\om =0$ so
that $uv\om =0$ (where $\om$ was defined in the Introduction). It
is easy to check that the last map induces a map $\psi
:D_{2k}V\otimes H^m (DV)\to H^{m+2k}(DV)$, $u\otimes (v+{\rm
Im}\,\partial_{m-2})\mapsto uv+{\rm Im}\,\partial_{m+2k-2}$. We
claim that $\psi (N\otimes H^m (DV))=0$, where $N$ is defined in
Lemma 3.3.

Indeed, let $u=x^{(a_1 )}_{1}y^{(b_1 )}_{1}\cdots x^{(a_m
)}_{m}y^{(b_m )}_{m}\in N$, where $\sum\limits^{m}_{i=1}(a_i +b_i
)=2k$. Then at least one exponent of $u$ is odd. For simplicity in
notation, let us assume that this exponent is $a_1$. Let
$v=x^{(c_1 )}_{1}y^{(d_1 )}_{1}\cdots x^{(c_m )}_{m} y^{(d_m
)}_{m}+{\rm Im}\,
\partial_{m-2}\in H^m (DV)$ be a basis element (Corollary 2.4), so
that $c_i +d_i =1$ for all $i$. Then
\begin{align}
    \psi(u\otimes(v+{\rm Im}\,\partial_{m-2}))= \notag \\
    =q x^{(a_1 +c_1 )}_{1}&y^{(b_1 +d_1 )}_{1}\cdots x^{(a_m +c_m
    )}_{m}y^{(b_m +d_m )}_{m}+ {\rm Im}\,\partial_{m+2k-2}, \tag{8}
\end{align}
where
\[
q=\binom{a_1 +c_1}{c_1} \binom{b_1 +d_1}{d_1}\cdots \binom{a_m
+c_m}{c_m}\binom{b_m +d_m}{d_m}.
\]
If $q$ is even, there is nothing o prove.
Let $q$ be odd. Then, since we have $c_i+d_i=1$, it follows that
for each $i$ at least one of $a_i+c_i+1, b_i+d_i+1$ is even and in
$K$ we have
\[
    (a_i+c_i+1)(b_i+d_i+1)=0.\tag{9}
\]
From the fact that $a_1$ and $q$ are odd we obtain that both
$a_1+c_1, b_1+d_1$ are odd. Let $z\in D_{m+2k-2}V$ be
\begin{align}
    z=x^{(a_1 +c_1-1 )}_{1}y^{(b_1 +d_1-1 )}_{1}x^{(a_2 +c_2 )}_{2}y^{(b_2 +d_2 )}_{2}\cdots x^{(a_m +c_m
    )}_{m}y^{(b_m +d_m )}_{m}.\notag
\end{align}
By a direct computation and using (9) we obtain
\begin{align}
    \partial_{m+2k-2}(z)=(a_1+c_1)(b_1+d_1)x^{(a_1 +c_1 )}_{1}y^{(b_1 +d_1 )}_{1}\cdots &x^{(a_m +c_m
    )}_{m}y^{(b_m +d_m )}_{m}\notag
    \\
    =x^{(a_1 +c_1 )}_{1}y^{(b_1 +d_1 )}_{1}\cdots x^{(a_m +c_m
    )}_{m}y^{(b_m +d_m )}_{m}.\tag{10}
\end{align}
From (8) and (10) we have that $\psi(u\otimes(v+{\rm
Im}\,\partial_{m-2}))=0$ and the claim is proved.

We thus have map of $G$-modules $D_{2k}/N \otimes H^m (DV)\to
H^{m+2k}(DV)$, such that $(u+N)\otimes(v+{\rm
Im}\,\partial_{m-2})\mapsto uv+{\rm Im}\,\partial_{m+2k-2}$. Using
Lemma 3.3 we obtain a map $ (D_k V)^{(1)} \otimes H^m (DV)\to
H^{m+2k}(DV)$. Using Corollary 2.4 one easily checks that this
last map is onto and hence an isomorphism because we have
\[
dim\big((D_kV)^{(1)}\otimes H^m(DV)\big)=\binom{2m+k-1}{k} 2^m
=dimH^{m+2k}(DV)
\]
by Corollary 2.5. Finally we recall that $D_k V=\Delta (k\om_1 )$
and, by Proposition 3.1, $H^m (DV)\simeq L(\om_m)$.
\end{proof}
\begin{corollary}\label{co3.4}
$H^{m+2}(DV)$ is an irreducible $G$-module of highest weight
$2\om_1 +\om_m$.
\end{corollary}
\begin{proof}
By Theorem 3.2, $H^{m+2}(DV)\simeq \De (\om_1 )^{(1)}\otimes
L(\om_m )$. But $\Delta (\om_1 )$ is just the natural module of
$G$, so that $\Delta (\om_1 )^{(1)}=L(\om_1 )^{(1)}$. By
Steinberg's tensor product theorem (see, for example, \cite[II.
3.17]{3}), $H^{m+2}(DV)\simeq L(\om_1 )^{(1)}\otimes L(\om_m
)\simeq L(2\om_1 +\om_m )$.
\end{proof}
\begin{rem} We have seen that the first two nonzero $H^{k}(DV)$ are
irreducible representations of $G$. It is not true in general that
all the nonzero $H^k (DV)$ are irreducible representations of $G$.
For example, let $m=1$ so that $G=SL(2,K)$. Then for each odd $k$
we have $H^k (DV)=D_k V$ from Lemma \ref{le2.1}. It is well known
that $D_k V$ is irreducible if and only if $k=2^n -1$, $n\ge 0$.
\end{rem}

\bigskip
\noindent Mihalis Maliakas\\
Department of Mathematics\\
University of Athens\\
Panepistimiopolis 15784\\
Athens, Greece\\
e-mail: mmaliak@math.uoa.gr
\end{document}